# Extraction de la racine carrée d'un entier naturel chez al-Baghdādī

Yassine Hachaïchi[1], Leïla Hamouda[2], Seïf Toumi[3]

## I. Introduction

Entre les neuvième et quinzième siècles, plusieurs mathématiciens arabes ont étudié des algorithmes numériques sur les entiers. Le premier d'entre eux, al-Khwārizmī (mort vers 850) les expose dans le cadre de la théorie du système décimal de position dans son ouvrage le calcul indien (*al-ḥisāb al-hindī*), ce livre est perdu dans sa version arabe et les versions latines ont été éditées par A. Allard[1].

L'extraction de la racine carrée d'un entier, qui est l'objet de cet article, repose sur un algorithme connu depuis al-Khwārizmī. Ce dernier, en effet, dans son livre « le calcul indien », dit que « *la manière de chercher la racine est connue* » [1, pp 52-56].

Bien qu'il se soit moins démocratisé que les quatre opérations usuelles de l'arithmétique, l'algorithme de l'extraction de la racine carrée a connu une postérité toujours actuelle, puisqu'il est implémenté aujourd'hui dans nos calculatrices et nos ordinateurs.

On retrouve cet algorithme dans les écrits de plusieurs mathématiciens comme Mūḥammad Ibn Mūsā al-Khwārizmī, dans son livre *Al-ḥisāb al-hindī* [1, pp 52-56], Abū'l-Ḥasan Aḥmad ibn Ibrāhim al-Uqlīdīsī (mort vers 980) dans son *Kitāb al-fūṣūl fī al-ḥisāb al-hindī* [2, pp 105-110], ʾAbū ʿAlī al-Ḥasan ibn al-Ḥasan ibn al-Haytham (mort vers 1040) dans un texte intitulé *Maqāla fī ʿillat al-jidhr wa idāfīhi wa naqlihi* [3, pp 463-467], ʾAbū'l-Ḥasan Kūshyār ibn Labbān ibn Bashahrī Gīlānī (mort vers 1030) dans *Kitāb fī ūsūl ḥisāb al-hind* [4], ʾAbū Mansūr ʿAbd al-Qāhir ibn Tāhir ibn Mūḥammad ibn ʿAbd-allāh al-Tamīmī al-Shāfiʿī al-Baghdādī (mort vers 1037) dans son traité *Al-takmila fī-l-ḥisāb* [5, pp70-82] et finalement Ghiyāth al-Dīn Jamshīd Masʿud al-Kāshī (mort vers 1430) dans son livre *Miftāḥ al-ḥisāb*[6].

Tous ces mathématiciens ont écrit sur l'algorithme dans le cas d'une racine entière d'un entier. Pour le cas de la racine non entière, nous trouvons des méthodes d'approximation dans les textes d'al-Khwārizmī [1, pp 52-56] et d'al-Baghdādī [5, pp70-82] ainsi que dans le texte d'Ibn al-Haytham. L'une est conventionnelle, traditionnellement utilisée par les calculateurs [3, pp 463-467] [5, pp70-82], l'autre est introduite par al-Khwārizmī [1, pp 52-56].

[1] **Université de Carthage, Ecole Nationale d'Ingénieurs de Carthage, LR18ES44, Laboratoire de Recherche Electricité Intelligente & TIC, EITIC.**
[2] [3] **Université de Tunis El Manar, Ecole Nationale d'Ingénieurs de Tunis, LR99ES20, Laboratoire de modélisation mathématique et numérique dans les sciences de l'ingénieur, équipe histoire des mathématiques, 2092, Tunis, Tunisie.**

Nous avons choisi d'analyser, dans cet article, le travail d'al-Baghdādī, sur l'extraction de la racine carrée des entiers, qui nous a semblé plus complet et plus détaillé que les travaux de ses prédécesseurs. A. Saïdane nous a fourni une courte biographie d'al-Baghdādī, il a également édité son ouvrage « *Al-takmila fī l-ḥisāb* » [5], où ce dernier expose dans un chapitre intitulé « *Fī bayāni kayfyati 'ikhrāji al-jūdhūri min al-a'dādi al-ṣiḥāḥi* » [5, pp 71-82] l'algorithme de l'extraction de la racine carrée des entiers. Ce chapitre est composé de six sections que nous transcrirons toutes. La transcription dans un langage mathématique moderne a pour but d'analyser l'algorithme. Il va de soi que cette transcription est porteuse de toute une *mathesis* anachronique, ceci nous engage ainsi que le lecteur, à faire la part entre ce qui appartient à al-Baghdādī et ce qui découle de notre transcription.

## II. Extraction de la racine carrée d'un entier naturel chez al-Baghdādī

Dans le chapitre « *Fī bayāni kayfyati 'ikhrāji al-jūdhūri min al-a'dādi al-ṣiḥāḥi* » [5, pp 71-82], al-Baghdādī présente l'extraction de la racine carrée d'un entier comme une méthode purement numérique. L'algorithme présenté est le même que dans les textes d'al-Khwārizmī et d'al-Uqlidīsī. Il donne des exemples bien choisis qui, en les comparant avec ceux pris par al-Khwārizmī, nécessitent un plus grand nombre d'itérations dans la recherche de leurs racines carrées, ce qui est dû à leurs grands nombres de décimaux. Dans chaque section, al-Baghdādī commence par donner une partie théorique et ensuite il valide sur des exemples. Vu la clarté de la manière avec laquelle al-Baghdādī présente son travail, il est probable que l'auteur veut que son exposé soit pédagogique. En effet, al-Baghdādī, après avoir quitté Baghdad, s'installa à Nichapour où il enseigna les mathématiques ainsi que les sciences religieuses [5, pp 10-11]. Nous pouvons, donc, penser que l'aspect pédagogique occupait une place importante dans son approche mathématique des problèmes à résoudre.

### 1. *Extraction de la racine carrée d'un carré parfait*

Cette première section du chapitre [5, pp 72-76] présente l'algorithme qui sera utilisé dans toutes les autres sections. D'une façon explicite, al-Baghdādī décrit les étapes à suivre pour extraire la racine d'un carré parfait. Sa description est itérative, dans le sens où il donne les deux premières itérations explicitement et il illustre la méthode sur trois exemples qui requièrent trois itérations. Vu qu'il utilise un *takht*, le résultat est obtenu à la dernière itération directement, sans avoir en vue les étapes antérieures. En effet, un *takht* est une « table à poussière » ou un bac à sable, les chiffres y sont tracés avec le doigt ou avec une baguette, les étapes de calcul étant facilement effacées au fur et à mesure que l'on avance dans le calcul.

Dans un premier temps, al-Baghdādī étiquette par « racine » les rangs impairs (en partant de l'unité) et par « non racine » les rangs pairs ; n'oublions pas que dans la représentation décimale, les rangs impairs correspondent aux puissances paires de 10. Il demande ensuite de se placer *sous* le dernier rang impair et de trouver la plus grande unité dont le carré est inférieur à ce qu'on a *au dessus*. Al-Baghdādī écrit :

> *Si tu veux extraire la racine d'un carré parfait, mets-le dans un* takht, *ensuite dis, à partir du premier rang* (de droite à gauche) *« racine », « non racine », « racine », etc… : aux rang impairs on aura « racine » et aux rang pairs, on aura « non racine »* [Les rangs impairs correspondent aux puissances paires et inversement, on commence à partir de $10^0$] *ensuite regarde la dernière position à laquelle tu as dit : « racine », et mets au-dessous le plus grand entier multiplié par lui-même qu'on peut soustraire de ce qu'il y a au-dessus. Si tu l'as connu alors mets-le au dessous, et multiplie-le par lui-même et soustrais-le de ce qu'il y a au-dessus, ensuite double ce nombre et décale son double d'un rang vers la droite, et regarde quel rang avant celui dans lequel tu as écrit alors met au-dessous de ce rang, le plus grand entier que tu multiplie par le double du rang déjà calculé, on peut le soustraire de ce qu'il y a au dessus, ensuite il est multiplié par lui-même on peut le soustraire de ce qu'il y a au-dessus, si tu l'as connu, mets-le au dessous, ensuite multiplie- le par le double, et ce que tu obtiens soustrait-le de ce qu'il y a au-dessus, ensuite multiplie-le par lui-même et soustrait-le de ce qu'il y a au-dessus. Puis double ce nombre aussi, et mets-le au-dessous à un rang décalé d'une case à droite. Ensuite, sur ce principe, on a la validation, la soustraction, le doublement et le transfert, jusqu'à ce que le nombre dont la racine est demandée se termine. Et si tu ne trouves pas à certains rangs de la validation et de l'écriture un nombre qu'on peut multiplier par celui qui le suit, et par lui-même, et qu'on peut soustraire de ce qu'il y a au dessus à chaque position alors mets un zéro sous cette position devant les multiples et décale la ligne au dessous d'un rang à droite et soumet toi au schéma déjà mentionné.*
>
> *Si le nombre dont la racine est demandée se termine, alors divise par deux toutes les positions que tu as doublées, alors ce qui reste, à la dernière ligne est la racine demandée. Si le nombre dont la racine est demandée se termine et il reste un nombre pair de zéros, alors partage en deux les positions que tu as doublées, et rajoute au début de la ligne* [de droite à gauche] *la moitié du nombre de zéros qui restent, alors ce qui reste est la racine demandée.*

Soit $N$ l'entier naturel dont on veut calculer la racine carrée.

Posons $N = n_k n_{k-1} \cdots n_1 n_0$ la décomposition décimale de $N$, ce qui équivaut à dire que

$$N = \sum_{i=0}^{k} 10^i \, n_i$$

avec $n_k, n_{k-1}$ non tous les deux nuls et les $n_i$ appartiennent à $\{0,1,2,\ldots,9\}$.

Si $k$ est pair alors il existe un entier $p$ tel que $k = 2p$.

Si $k$ est impair alors on rajoute un 0 à gauche de $N$ et on se ramène au cas où $k$ est pair. Al-Baghdādī ne mentionne pas explicitement cela mais en demandant à se placer au-dessous du rang impair dans son étiquetage racine-non racine, nous déduisons que cela peut être effectué implicitement.

Dans la suite, pour expliquer l'algorithme, nous traiterons le cas où $k$ est pair.

Pour un carré parfait, on note $s$ sa racine carrée écrite sous la forme décimale $s = \sum_{i=0}^{p} 10^i \, s_i$

Nous noterons, également, le nombre $s$ dont les décimales sont $s_p, s_{p-1},\ldots, s_0$ par

$$(s_p, s_{p-1}, \ldots, s_1, s_0)$$

D'abord, al-Baghdādī cherche $s_p$ le plus grande décimale telle que :

$$s_p^2 \leq (n_k, n_{k-1}) < (s_p + 1)^2.$$

Il soustrait alors $s_p^2$ de $(n_k, n_{k-1})$ et décale d'un rang vers la droite, avec nos notations au dessous de $n_{k-2}$, et pose $2\, s_p$.

D'ailleurs en reprenant l'exemple du manuscrit $N = 54756$, al-Baghdādī commence par étiqueter : **$0^n 5^r\ 4^n 7^r 5^n 6^r$** en notant par $n$ non racine et $r$ par racine. Le **$0^n$** est rajouté dans notre travail pour mieux expliquer l'algorithme, en effet, n'importe quel nombre peut être considéré comme une suite de couples (racine, non racine). Pour cet exemple, on a $k = 6$ donc $p = 3$. De plus, en remarquant que $4 \leq 5 < 9$, il prend $s_3 = 2$. Il double alors le 2 qu'il a mis et il le décale d'un rang à droite, il obtient :

| 1 | 4 | 7 | 5 | 6 |
|---|---|---|---|---|
|   | 4 |   |   |   |

On notera également la multiplication par « · ».

Ensuite, il cherche le plus grand nombre $s_{p-1}$ tel que $s_{p-1} \cdot (2.s_p, s_{p-1})$ soit inférieur à $(n'_{2p}, n'_{2p-1}, n_{2p-2})$ où le nombre $(n'_{2p}, n'_{2p-1})$ est $(n_{2p}, n_{2p-1})$ auquel il a retranché $s_p^2$.

Dans l'exemple précédent, on a $3 \cdot (43) = 129 \leq 147$ et $4 \cdot (44) = 176 > 147$, on a alors $s_2 = 3$. Il retranche ensuite $s_{p-1} \cdot (2.s_p, s_{p-1})$ de $(n'_{2p}, n'_{2p-1}, n_{2p-2})$, il double $s_{p-1}$ et il décale le tout d'un rang à droite, ce qui nous donne dans l'exemple :

| 1 | 8 | 5 | 6 |
|---|---|---|---|

|   | 4 | 6 |   |
|---|---|---|---|

Il continue ce procédé tel que $s_{p-i} \cdot (2.s_p, 2.s_{p-1}, \cdots, s_{p-i})$ soit inférieur à $(n'_{2p}, n'_{2p-1}, \cdots, n'_{2p-2i})$, il retranche le premier terme cité du second, il double $s_{p-i}$ et il décale le tout d'un rang à droite, jusqu'à obtenir 0. Ce qui nous donne dans l'exemple :

$4 \cdot (464) = 1856$, en le retranchant de 1856, il reste 0.

| 0 | 0 | 0 | 0 |
|---|---|---|---|
|   | 4 | 6 | 4 |

Pour obtenir la racine, al-Baghdādī divise le nombre obtenu par deux, en dehors de l'unité, et il obtient la racine qui est 234 dans son exemple. On résume l'exécution de l'exemple précédent dans le tableau suivant :

| 5 | 4 | 7 | 5 | 6 |
|---|---|---|---|---|
| 2 |   |   |   |   |
| 1 | 4 | 7 | 5 | 6 |
|   | 4 |   |   |   |
|   |   | 3 |   |   |
|   | 1 | 8 | 5 | 6 |
|   |   | 4 | 6 |   |
|   |   |   |   | 4 |
| 0 | 0 | 0 | 0 | 0 |
|   |   | 4 | 6 | 4 |

Dans l'exemple suivant, il considère le cas où $s_{p-i} \cdot (2.s_p, 2.s_{p-1}, \ldots, s_{p-i})$ est supérieur à $(n'_{2p}, n'_{2p-1}, \cdots, n'_{2p-2i})$ pour toute valeur de $s_{p-i} \geq 1$, alors on prend $s_{p-i} = 0$. Voici l'exemple du texte :

| 4 | 1 | 2 | 0 | 9 |
|---|---|---|---|---|
| 2 |   |   |   |   |
| 0 | 1 | 2 | 0 | 9 |
|   | 4 |   |   |   |
|   |   | 4 | 0 |   |
|   |   |   |   | 3 |
| 0 | 0 | 0 | 0 | 0 |
|   |   | 4 | 0 | 3 |

La racine obtenue est alors 203.

A la fin de cette section, il considère le cas où il reste un nombre pair de zéros et que le procédé n'est pas encore terminé, il suggère alors de diviser le nombre de zéros restants par deux et de le mettre à la fin de la valeur obtenue. Il prend pour exemple la racine de 5290000

| 5 | 2 | 9 | 0 | 0 | 0 | 0 |
|---|---|---|---|---|---|---|
| 2 |   |   |   |   |   |   |
| 1 | 2 | 9 | 0 | 0 | 0 | 0 |
|   | 4 |   |   |   |   |   |
|   | 3 |   |   |   |   |   |
|   | 0 | 0 | 0 | 0 | 0 |   |
|   | 4 | 3 |   |   |   |   |

Il prend alors la racine 23 à laquelle il rajoute deux zéros et obtient 2300, comme racine de 5290000.

## 2. *Extraction de la racine carrée d'un carré non parfait*

Dans cette partie, al-Baghdādī écrit :

*Extraire sur le* takht *la racine d'un non carré est comme extraire de la racine d'un carré sauf qu'il reste dans le carré non parfait des fractions, une partie ou des parties après l'extraction de la racine de la partie entière. Les calculateurs ne sont pas d'accord sur l'attribution des parties restantes : Mūḥammad Ibn Mūsā al-Khwārizmī, qu'il ait la miséricorde de Dieu, l'attribue au double de la racine sortante des entiers et la plupart des calculateurs l'attribuent au double de la racine sortante plus un. Et ce dire est le plus juste.*

*Exemple : si nous voulons la racine de cent cinquante cinq, sous cette forme 155, alors on extrait sa racine avec approximation, en suivant le dessin précédent sur l'extraction des racines, alors il reste du calcul sous cette forme 11*

*22*

*après on divise par deux le deux qui est à la deuxième position de la ligne au-dessous, car il était un qu'on a doublé, alors il reste dans la ligne du dessous douze, ensuite on met ce douze au-dessus du onze des fractions, après on double le douze, on rajoute à son double un, à la fin on lui attribue les parties restantes alors on sait que la racine du nombre qu'on a pris dans cet exemple est douze dirhams et onze portions des vingt cinq portions du dirham sous cette forme 12*

*11*

*25*

Dans ce qui suit, al-Baghdādī donne l'exemple de l'extraction de la racine carrée de 2 comme contre-exemple à l'approximation d'al-Khwārizmī ; on montre dans l'article [10] que ce contre-exemple est erroné.

*Et ce qui prouve que les calculateurs ont raison d'attribuer les parties restantes de la racine d'un nombre non carré au double de la partie entière de la racine extraite auquel on a ajouté un un, la non véracité de l'affirmation de Mūḥammad Ibn Mūsā al-Khwārizmī quand il attribue les*

*parties au double de la partie entière de la racine extraite seulement est : si nous voulons extraire la racine de deux par approximation, on obtient un et il reste de ses parties un, alors si on rapporte au double de un, comme dit Ibn Mūsā, alors la racine serait un-et-demi ; alors que si on le rapporte au double de la racine sortante, auquel on ajoute un un, comme le disent la majorité, la racine devient un-et-un-tiers. Si on multiplie le un-et-demi par lui-même alors on obtient deux-et-un-quart et si on multiplie le un-et-un-tiers par lui-même, on obtient un-et-sept-neuvième. Et entre ce dernier et deux, deux neuvièmes manquent pour compléter deux et ceci est plus proche de deux que du quart qui est en excès à deux.*

Al-Baghdādī prend également l'exemple de 3 qui est lui aussi erroné et dit « alors on sait par cela la véracité du dire de la majorité des calculateurs ».

Après transcription, on note l'entier $N$ un carré non parfait c'est-à-dire que

$$N = E^2 + R \qquad (avec\ 0 \leq R \leq 2E)$$

où

$$E = \sum_{i=0}^{n} s_i\, 10^i \qquad et \qquad R = N - \left(\sum_{i=0}^{n}(s_i 10^i)^2 + 2\sum_{i,j=0, i>j}^{n} s_i 10^i s_j\, 10^j\right)$$

Il cherche $E$ par le même algorithme que la première partie. Finalement, il estime la racine $s$ par

$$s \cong E + r \qquad où \qquad r = \begin{Bmatrix} \dfrac{R}{2E+1} \\ ou\ bien \\ \dfrac{R}{2E} \end{Bmatrix}$$

Si $r$ est estimé à $\dfrac{R}{2E}$ alors c'est l'approximation proposée par al-Khwārizmī [1, pp 52-53].

Si $r$ est estimé à $\dfrac{R}{2E+1}$ alors c'est l'approximation dite conventionnelle considérée par les calculateurs [5, p76].

Al-Baghdādī expose les deux approximations et conclut hâtivement que l'approximation conventionnelle du reste est meilleure en donnant deux exemples d'extraction de racine des nombres 2 et 3.

Pour son premier exemple, al-Baghdādī écrit 2=$1^2$+1 et considère les deux approximations $\sqrt{2} \approx 1 + \dfrac{1}{2}$ et $\sqrt{2} \approx 1 + \dfrac{1}{3}$ puis il les élève au carré pour avoir successivement : 2+$\dfrac{1}{4}$ et 1+$\dfrac{7}{9}$. Il constate que $1 + \dfrac{7}{9}$ est plus proche de 2 que 2+$\dfrac{1}{4}$. Il considère ensuite 3=$1^2$+2 et considère les deux approximations $\sqrt{3} \approx 2$ et $\sqrt{2} \approx 1 + \frac{2}{3}$ puis il les élève au carré pour avoir successivement : 4 et $2 + \frac{7}{9}$. Il constate que $2 + \frac{7}{9}$ est plus proche de 3 que 4. Il en conclut que ceci reste vrai pour tout entier.

Pour montrer que cette conclusion est fausse dans le cas général, prenons l'exemple de $10=3^2+1$. Par la méthode d'al-Khwārizmī, on obtient $\sqrt{10} \approx 3 + \frac{1}{6}$ et par celle dite conventionnelle $\sqrt{10} \approx 3 + \frac{1}{7}$, en les élevant respectivement au carrée, on obtient $10 + \frac{1}{36}$ pour la première et $9 + \frac{43}{49}$ pour la seconde. Or, $0 < \frac{1}{36} < \frac{6}{49}$, qui nous donne une meilleure approximation par la première méthode.

Nous remarquons, donc, que le choix d'approximation d' al-Baghdādī pour un entier naturel $N$ qui n'est pas un carré parfait, n'est pas bien précis. En effet, il conclut, par une justification empirique, que l'approximation conventionnelle est meilleure que celle d'al Khwārizmī alors que nous montrons dans notre article [9] « Extraction de la racine carrée d'un entier chez Ibn al Haytham et comparaison avec al-Baghdādī» que deux cas se présentent. Si $R \leq E - 1$ alors l'approximation d'al Khwārizmī est meilleure pour l'estimation de la racine carrée et si $R \geq E$ alors l'approximation conventionnelle est meilleure.

Pour les deux exemples exposés par al-Baghdādī, $N = 2$ et $N = 3$, on remarque que dans ces deux cas, $R \geq E$. En effet, dans le cas où $N = 2$, $E = 1$ et $R = 1$ et dans le cas où $N = 3$, $E = 1$ et $R = 2$.

### 3. *Extraction de la racine d'un nombre en le multipliant par un nombre*

La troisième partie introduit essentiellement la technique d'extraction de la racine carrée d'un entier $N$ en le multipliant par un autre entier $A$ un nombre pair de fois de telle manière à avoir

$$\sqrt{A^{2p}N} = A^p\sqrt{N}$$

Al-Baghdādī commence par valider sa méthode de calcul sur des exemples de carrés parfaits en prenant l'exemple suivant:

$4 \times (3 \times 3 \times 3 \times 3) = 324$ dont la racine est 18, en divisant 18 par ($3 \times 3 = 9$, $p = 2$) on obtient 2 qui est la racine de 4.

Il reprend un autre exemple dans le même esprit :

$4 \times (3 \times 3) \times (5 \times 5) = 900$ dont la racine est 30, en divisant 30 par ($3 \times 5 = 15$, $p = 1$) on obtient 2 qui est la racine de 4.

Il considère ensuite un nombre dont la racine n'est pas entière : 2, ensuite il le multiplie par 3 quatre fois ($2 \times 3 \times 3 \times 3 \times 3 = 162$), il applique l'algorithme d'extraction de la racine entière et obtient :

$162 = 12^2 + 18$. Si on prend la racine entière, à savoir 12, et on la divise par (3×3=9) on obtient $\sqrt{2} \approx 1 + \frac{1}{3}$.

4. ***A propos de l'extraction de la racine carrée en multipliant par des puissances de 10***

Cette méthode consiste à chercher la racine carrée d'un entier *N* en le multipliant par un nombre pair de 10, c'est-à-dire on prend $A = 10$ dans l'équation

$$\sqrt{A^{2p}N} = 10^p\sqrt{N} .$$

Les $10^p$ représentent les *p* chiffres après la virgule de la racine carrée. Al-Baghdādī traduit ces chiffres après la virgule en minutes et secondes (base 60) en multipliant par soixante à chaque fois. Il insiste sur le fait que plus *p* est grand, plus l'approximation de la racine carrée est précise.

L'exemple explicité par al-Baghdādī dans cette section est le calcul de la racine de 5, il rajoute six zéros, il obtient alors 5 000 000. En utilisant la méthode de la section 2, il obtient 5000000=2236²+304. Ensuite il décale de trois rangs (qui est la moitié des six zéros rajoutés), il prend le 2 du début qui est la partie entière, il reste 236, qu'il multiplie par 60 et il obtient 14160. Ensuite il décale de trois rangs, il prend 14 qui est le nombre de minutes, il reste 160 qu'il multiplie par 60 pour obtenir 9600. Il décale encore de trois rangs, il prend 9 qui est le nombre de secondes, il reste 600 qu'il multiplie par 60 pour obtenir 36000, qui équivaut 36 tierces.

2236=2×1000+236,

236×60=14×1000+160,

160×60=9×1000+600,

600×60=36×1000+0.

$\sqrt{5} = 2 + 14$ minutes + 9 secondes + 36 tierces et $\sqrt{5} = 2.236$ en décimal.

En outre, l'erreur de choix d'approximation réalisée par al-Baghdādī dans la section 2 est atténuée par le fait que dans cette section, en multipliant par une puissance paire de dix l'entier dont on cherche la racine carrée, il assure que le résultat qu'il obtient est très précis. En effet, le nombre de zéros qu'il ajoute se traduit, de nos jours, par le nombre de chiffres après la virgule dans le calcul de la racine carrée.

5. ***Vérification de l'exactitude du calcul***

La cinquième partie est une vérification de l'exactitude du calcul de la racine carrée modulo 9. Soit *N* le nombre dont on cherche la racine carrée et *s* sa racine carrée.

<u>Premier cas</u> : *N* est un carré parfait c'est-à-dire que $N = s^2$ où *s* est un entier. Après avoir cherché *s* par l'algorithme donné à la première section, on note *a* et *b* les restes respectifs de la division euclidienne de *N* et $s^2$ par 9. Si $a = b$ alors, d'après al-Baghdādī, le calcul de la racine carrée est juste.

Il prend l'exemple où $N = 3721$.

On a $3721 \equiv 4[9]$ *et* $3721 = 61^2$ or $61 \equiv 7[9] \Rightarrow 61^2 \equiv 4[9]$.

3721 et $61^2$ ont 4 comme reste dans la division euclidienne par 9 donc le calcul de la racine est juste.

<u>Deuxième cas</u> : *N* n'est pas un carré parfait c'est-à-dire que

$$N = E^2 + R \quad (avec\ 0 \leq R \leq 2E)$$

et sa racine *s* s'écrit comme $s \cong E + r$ où $r = \left\{ \begin{array}{c} \frac{R}{2E+1} \\ ou\ bien \\ \frac{R}{2E} \end{array} \right\}$

Soient *a*, *b* et *c* les restes respectifs de la division euclidienne des entiers *N*, $E^2$ et *R* par 9 c'est-à-dire que $N \equiv a[9]$, $E^2 \equiv b[9]$ *et* $R \equiv c[9]$.

Alors al-Baghdād affirme que si $a = b + c$ alors le calcul de la racine carrée est juste.

Il considère l'exemple de $N = 249$, il obtient une racine entière de 15 et comme reste 24. $249 \equiv 6[9]$ or $15^2 \equiv 6^2[9] \equiv 0[9]$ et $24 \equiv 6[9]$ d'où $15^2 + 24 \equiv 6[9]$.

Dans cette partie, al-Baghdādī prend une implication, la non-conformité modulo 9 implique une erreur de calcul, pour une équivalence. Il écrit [5, pp.80-81]:

> *S'il lui correspond alors le travail est juste, s'il ne lui correspond pas alors c'est faux.*

### 6. *A propos des critères de vérification de l'extraction des racines*

La dernière partie étudie les décimales qui composent les carrés parfaits (modulo 10 et modulo 9). Dans cette partie, al-Baghdādī affirme certaines propriétés des racines carrées extraites des carrés parfaits.

L'unité d'un carré parfait ne peut être que 0, 1, 4, 5, 6 ou 9.

Si l'unité d'un carré parfait est 0 alors sa racine aussi a 0 pour unité, si c'est un 5 alors sa racine aussi a 5 pour unité, si c'est un 1 alors l'unité de la racine est soit 1 soit 9, si c'est un 4 alors l'unité de la racine est soit 2 soit 8, si c'est un 6 alors l'unité de la racine est soit 4 soit 6, et finalement, si c'est un 9 alors l'unité de la racine est soit 3 soit 7.

Les carrés parfaits sont tous congrus à 0, 1, 4, 7 modulo 9. Al-Baghdādī utilise pour l'illustrer la méthode de la somme des décimales.

### 7. *Efficacité de l'algorithme*

#### a. *Utilisation actuelle de l'algorithme*

Le travail d'al-Baghdādī montre que l'algorithme de l'extraction de la racine carrée est déjà bien assimilé et étudié, du point de vue numérique, au onzième siècle déjà. Son efficacité et sa « simplicité » ont contribué à la fortune de cet algorithme dont la large diffusion et son actuelle utilisation sont témoins. En effet, cet algorithme est encore utilisé de nos jours et des améliorations lui ont été apportées ainsi, dans les articles [7] et [8], les auteurs améliorent la version *hardware* de cette méthode d'extraction de la racine carrée, en optimisant l'espace mémoire utilisé dans chaque article. Le *takht*, utilisé par al-Baghdādī, et l'implémentation binaire, utilisée de nos jours, constituent deux facettes différentes d'un même objectif : utiliser le minimum d'espace matériel pour le premier et de place mémoire pour le second.

#### b. *Algorithme d'extraction de la racine carrée et la méthode de Newton*

La version *software* d'extraction de la racine carrée, celle de Newton, est basée sur la méthode de point fixe : $u_0 = a$ et $u_{n+1} = \frac{u_n^2 + a}{2\,u_n}$, dont la convergence est rapide. En choisissant $a$ un entier et en comparant avec la méthode d'al-Baghdādī, nous avons trouvé que cette dernière fournit une meilleure précision que celle de Newton qui converge après quatre itérations (pour les entiers entre 0 et 1023).

### III. Conclusion :

L'originalité du travail d'al-Baghdādī demeure dans le fait d'être, non seulement, pédagogique et très explicatif mais en plus, purement numérique. Son approche va dans le même sens que celui initié par al Khwārizmī, de traitement syntaxique d'un problème mathématique. Comme un linguiste, al-Baghdādī s'assure que l'ensemble des règles des énoncés est bien formé – il procède en cela en algébriste- et qu'ils ont bien un sens. En revanche, il ne s'intéresse pas à la justification de l'algorithme qu'il utilise ni à sa bien-fondé mathématique.